\documentclass[11pt]{article}

\usepackage{amssymb}
\usepackage{amsfonts, eufrak}
\usepackage[dvips]{graphics}
\usepackage[latin1]{inputenc}
\usepackage{amsmath}

\def\dd{\displaystyle}

\def\ll{\lambda}
\def\dd{\displaystyle}

\def\rr{\mathbb{R}}

\newtheorem{thm}{Theorem}

\begin{document}
\title{An elementary approach to certain bilinear estimates}
\author
    {Jos\' e A. Barrionuevo%
     \thanks{josea@mat.ufrgs.br} \\
     Lucas Oliveira%
     \thanks{lucas.oliveira@ufrgs.br} \\
     Departamento de Matem\' atica \\
     UFRGS \\
     Av. Bento Gon\c calves 9500, 91509-900  Porto Alegre, RS, Brasil \\
     Jarod Hart%
     \thanks{jhart@math.ku.edu} \\
     Department of Mathematics \\
     University of Kansas \\
     Lawrence, Kansas 66045-7594, USA}

\maketitle
\begin{abstract} We prove $L^2(\rr )\times L^2(\rr ) \to L^{1,\infty}(\rr )$ estimates for some bilinear maximal operators of Kakeya and lacunary type. Our method is geometric and elementary and can possibly be applied to other situations. This is a first draft, so comments are welcome.
\end{abstract}

\newpage

\section{Introduction}
Consider the following operator, initially defined for Schwartz functions $f(x), g(x) \in \mathcal{S}(\rr )$ and $\delta > 0$ small,
\begin{equation}
\mathfrak{M}_{\delta}(f,g)\, (x) = \sup_{R\in \mathcal{B}_x} \frac{1}{| R |}\, \int_{R} \, | F(y,z) |\, dy dz
\end{equation}
where $\dd F(y,z) = f(y)\, g(z) \in \mathcal{S}(\rr^2 )$, $\dd \mathcal{B}_x$ is the class of all $\dd 1 \times \delta$ rectangles in $\rr^2$ centered at $(x,x)$ having longest side pointing along a $\delta$-separated set of directions $\Omega$. Let $\dd D= \{ (x,x) : x \in \rr \} \subset \rr^2$ be endowed with the standard $1-$dimensional Lebesgue measure. We will use the following notation. If $E$ is a subset of $\rr$ or $D$, then $\dd |E|$ denotes its Lebesgue measure while for a rectangle in $\rr^2$ it is the $2$-dimensional Lebesgue measure.
\begin{thm} We have the following $\dd L^2 \times L^2 \to L^{1,\infty}$ estimate: for $\ll > 0$, let \\
$\dd E_{\ll} = \{ x \in \rr:\, |\, \mathfrak{M}_{\delta}(f,g)(x) | \geq \ll \}$. Then
\begin{equation}
|E_{\ll} | \lesssim \left( \log\left(\frac{1}{\delta}\right)\right)^{1/2}\, \frac{1}{\ll}\; ||f||_2\, ||g||_2
\end{equation}
\end{thm}

{\bf Remark.} Using the trivial pointwise estimate $\dd \mathfrak{M}_{\delta}(f,g)(x) \leq \delta^{-1}\, M(f,g)(x)$, where $M$ is the Hardy-Littlewood type bilinear operator studied in \cite{BHO1}, which is valid for positive functions, one only gets
\[ |E_{\ll} | \lesssim \frac{1}{\delta^{1/2}}\, \frac{1}{\ll}\; ||f||_2\, ||g||_2 \]
We do not know if the estimate $(2)$ is sharp.

The method of proof allows us to extend to this bilinear setting, the result of Nagel, Stein and Wainger on lacunary
maximal operators.

\begin{thm}
Let $\dd M_{Lac} (f,g) (x)$ be as in (1) but with $\dd\mathcal{B}_x$ denoting the class of all rectangles in $\rr^2$
with longest making an angle of $2^{-j}$ with $D$. Then there is the $\dd L^2 \times L^2 \to L^{1,\,\infty}$ estimate
\begin{equation}
|\,\{ x \in D:\, M_{Lac} (f,g)(x) > \ll\,\}\, | \lesssim \ll^{-1}\,\, || f ||_2 || g ||_2
\end{equation}
\end{thm}

\section{Proofs}
We can assume that $f, g$ are positive and supported on $[-3,3]$ due to the local nature of $\mathfrak{M}_{\delta}$, which will, in turn, be supported on $[-5,5]$.

Now, given $\ll > 0$ and $x \in E_{\ll}$, there exists $R_x \in\mathcal{B}_x$ such that
\[ \frac{1}{| R_x |}\, \int_{R_x} \, F(y,z)\, dy dz > \ll \]
Define $\dd I_x = R_x \cap D$. We will identify $D$ with $\rr$ and, considering arbitrary compact subsets $K$ of $E_{\ll}$ we have that $K$ is covered by a finite family $I_j$, for $\dd j\in L$.

Applying Vitali's lemma we can select a disjoint sub family, which we still call $I_j$ such that
$\dd |\bigcup_j\, I_j | \geq \frac{1}{10}\, | E_{\ll} |$. From $\delta$-separation and elementary geometric considerations we have that \\
$\dd | I_j | \approx \left( 1 - 2^{s_j}\, | R_j |\,\right) = \left( 1 - 2^{s_j}\, \delta\right) =
(\delta^{-1} - 2^j)\, |R_j|$, where $\dd s_j \in \left\{ 0, 1,\cdots , \log\left(\frac{1}{\delta}\right)
\right\}$.

Let $\dd L_1 = \{ j\in L: \left( \delta^{-1} - 2^{s_j}\right)\, \leq 10 \}$, $L_2 = L - L_1$.
 Define $\dd K_i = K \cap\, \left( \cup_{j\in L_i} R_j \right) , i = 1,2$.

We estimate $\dd | K_1 |$
\begin{eqnarray}
|K_1| \lesssim \sum_{j\in L_1} | I_j | & = & \sum_{j\in L_1} (\delta^{-1} - 2^j)\, |R_j| \\
                               & \leq & 10\,\sum_{j\in L_1} \,\ll^{-1}\, \int_{R_j}F(y,z)\, dydz \\
                               & = & 10\,\, \ll^{-1} \,\int_{\rr^2}\,\left(\sum_{j\in L_1}  \chi_{R_j} (y,z)\,\right)\, F(y,z)\, dydz \\
                               & \leq & 10\,\, \ll^{-1} \, \left|\left| \sum_{j\in L} \, \chi_{R_j} (y,z) \right|\right|_2 \; \left|\left|\, F(y,z)\, \right|\right|_2
\end{eqnarray}
Observe that  $\# L = \mathcal{O}(\, \delta^{-1})$ and that for all $p \geq 1$
\begin{equation}
|| F ||_p = || F ||_{L^p(\rr^2)} = || f ||_{L^p(\rr )}\, ||g ||_{L^p(\rr )} = || f ||_p\, || g ||_p
\end{equation}
 Making use of C\' ordoba's estimate below
\begin{equation}
||\, \sum_{j\in L_1} \chi_{R_j} \, ||_2 \lesssim  \log\left(\frac{1}{\delta}\right)^{1/2}\,
\left( \sum_{j\in L_1} |R_j| \right)^{1/2}
\end{equation}
we obtain
\begin{equation}
|K_1| \lesssim  \ll^{-1} \, \log\left(\frac{1}{\delta}\right)^{1/2} \; || f ||_ 2\, ||g ||_2
\end{equation}

Before we estimate $|K_2|$, we need some simple results on an auxiliary maximal operator. For $l, w > 0$, $x\in \rr$, let $P_{x,l,w}$ be the parallelogram in $\rr^2$ with center $(x,x)$, with two vertical sides of length $2 w$ and the other two parallel to $D$ and length $l$ so that its vertices are $\dd (x+l, x+l+w),\, \\
(x+l, x+l-w),\, (x-l, x-l+w)$,
and $(x-l,x-l-w)$. Consider the maximal operator defined by
\begin{equation}
M_D(f,g)(x) = \sup_{h,w}\,\frac{1}{| P_{x,l,w} |}\, \int_{P_{x,l,w}} | F(y,z) |\, dydz
\end{equation}
If $M_1$ is the $1-$dimensional Hardy Littlewood operator and $M_V$ denotes the operator in $\rr^2$ acting on the vertical variable $z$ only, given by
\begin{equation}
M_V F(y,z) = \sup_w\, \frac{1}{2 w} \int_{-w}^{w} F(y,z+s)\, ds
\end{equation}
we have, observing that for $\dd f,g \geq 0$, $\dd M_V F (x,x) \lesssim f(x)\, M_1\, g(x)$, that
\begin{equation}
M_D(f,g)(x) \lesssim M_1 (f\, M_1\, g)(x)
\end{equation}
The above inequality together with the  classical $L^1\to L^{1,\infty}$ mapping property of $M_1$ implies the following $L^2 \times L^2 \to L^{1,\infty }$ estimate \begin{equation}
|\, \{ x\in\rr :\, M_D(f, g)(x) > \ll \}\, | \lesssim \ll^{-1}\, ||\, f\, M_1\, g\, ||_1 \lesssim
\ll^{-1}\; || f ||_ 2\, || g ||_2
\end{equation}

We are now ready to estimate $\dd | K_2 |$. The key observation is the simple geometric fact that the rectangles $R_j$ for $j\in L_2$ form an angle with $D$ bounded by $\dd C\, \delta$ where $C$ is an absolute constant ($C < 1000$).
This implies that we can find a parallelogram $P_j$ as above with vertical sides $\dd\approx C\,\delta$ and with the sides parallel to $D$ with length $\dd\approx 1$ such that $\dd R_j \subset P_j$ and that $\dd | P_j | / | R_j | \leq C_1$ where $C_1$ is another absolute constant.
\begin{equation}
\frac{1}{| R_j |}\,\int_{R_j} | F(y,z) |\; dy dz \leq C_1 \, \frac{1}{| P_j |}\,\int_{P_j} | F(y,z) |\; dy dz \end{equation}
But this implies that for another absolute constant $c > 0$ we have
\begin{equation}
 \bigcup_{j\in L_2} I_j\, \subset \{ x\in D: M_D (f, g)(x) > c \ll\,\}
\end{equation}
By $(14)$ this gives
\begin{equation}
| K_2 | \lesssim \ll^{-1}\,\, || f ||_2|| g ||_2
\end{equation}
Estimates (10) and (17) imply (2) proving Theorem 1.

The proof of Theorem 2 is even simpler. Given $\ll > 0$ and $ x\in E_{\ll}$, we obtain $R_x$, with
\begin{equation}
\frac{1}{| R_x |} \int_{R_x} F(y,z)\, dy dz > \ll
\end{equation}
If $K$ is any compact subset of $E_{\ll}$, it is covered by a finite collection $\dd \{ I_j \}$. We apply Vitali's lemma to get a disjoint sub family with measure $\geq c | E_{\ll}|$. We split the collection into two classes $L_1, L_2$ as before according to whether $| I_j | \leq 10 | R_j |$ or not. In the first case we repeat the calculation with the difference of using the linear lacunary estimate instead of Kakeya. The result is
\begin{eqnarray}
\sum_{j\in L_1} | I_j | & \leq & \sum_j 10\; \ll^{-1}\, \int_{R_j} F(y,z)\, dy dz \\
                        & \lesssim & \ll^{-1} \int_{\rr^2}
                        \left(\sum_{j\in L_1}\,\chi_{R_j}(y,z)\right) \; F(y,z)\; dy dz \\
                        & \leq & \ll^{-1} \; || \sum_{j\in L_1}\,\chi_{R_j}(y,z) ||_2\; || f ||_2\, || g ||_2 \\
                        & \lesssim & \ll^{-1} \; || f ||_2\, || g ||_2
\end{eqnarray}
The remaining $I_j$ can be controlled using the operator $M_D$ in a way similar to (15)-(17) and we leave the details to the interested reader.

\end{document}